\documentclass{amsart}
\usepackage{latexsym,amsmath,amsfonts,amssymb,amsxtra}
\usepackage{color,latexsym,amsmath,amsthm,amsfonts,amssymb,xcolor,verbatim,lipsum,layout}

\usepackage{enumerate}

\usepackage{soul,cite,lineno,cancel}
\modulolinenumbers[2]

\scrollmode

\newtheorem{thm}{Theorem}[section]
\newtheorem{deff}[thm]{Definition}
\newtheorem{lem}[thm]{Lemma}

\newtheorem{rem}[thm]{Remark}
\newtheorem{prop}[thm]{Proposition}

\newcommand{\ve}{\varepsilon}
\newcommand{\wt}{\widetilde}
\newcommand{\wh}{\widehat}

\def\vt{\vartheta}
\newcommand{\ov}{\overline}
\newcommand{\mg}{\marginpar}

\def\N{{\mathbb N}}

\def\S{\mathbb S}
\def\mcL{{\mathcal L}}

\def\mcP{{\mathcal P}}

\def\mcK{{\mathcal K}}
\def\mcM{\mathcal M}
\def\mcN{\mathcal N}
\def\vL{{\varLambda}}
\def\mfA{{\mathfrak A}}
\def\mfB{{\mathfrak B}}
\def\mfC{{\mathfrak C}}
\def\mfD{{\mathfrak D}}

\def\mfM{\mathfrak M}
\def\mfN{\mathfrak N}
\def\mfZ{\mathfrak Z}

\def\mfN{\mathfrak N}

\def\emp{\emptyset}

\def\mcN{\mathcal N}
\def\smi{\setminus}

\def\cd{\circledast}
\def\dd{\divideontimes}
\def\tr{\triangle}

\begin{document}

\title{Splitting of  Liftings in Product Spaces II}
\author[K. Musia{\l}]{Kazimierz Musia\l}
\address{50-384 Wroc{\l}aw, Pl. Grunwaldzki 2/4, Poland}
\email{kazimierz.musial@math.uni.wroc.pl}
\thanks{}

\date{\today}
\begin{abstract}
Let $(X, \mfA,P)$ and $(Y, \mfB,Q)$ be two probability spaces, $R$ be their skew product on the product $\sigma$-algebra $\mfA\otimes\mfB$ and  $\{(\mfA_y,S_y)\colon y\in{Y}\}$ be a $Q$-disintegration of $R$. Then let $\mfA\dd\mfB$ be the $\sigma$-algebra generated $\mfA\otimes\mfB$ and by the family $\mcM:=\{E\subset{X\times{Y}}\colon \exists\;N\in\mfB_0\;\forall\;y\notin{N}\;\wh{S_y}(E^y)=0\}$ and $\wh{R_{\dd}}$ be the extension of  $R$ such that  $\mcM$ becomes the family of $\wh{R_*}$-zero sets ($\wh{S_y}$ is the completion of $S_y$ and $\mfB_0=\{B\in\mfB: Q(B)=0\}$). We prove that  there exist a lifting $\pi$ on $\mcL^{\infty}(\wh{R_{\dd}})$  and liftings $\sigma_y$ on $\mcL^{\infty}(\wh{S_y})$ , $y\in Y$, such that
\[
[\pi(f)]^y= \sigma_y\Bigl([\pi(f)]^y\Bigr) \qquad\mbox{for every}
\quad y\in Y\quad\mbox{and every}\quad f\in\mcL^{\infty}(\wh{R_{\dd}}).
\]
In case of a separable $P$ and in case when $R\ll{P}\times{Q}$  a characterization of stochastic processes possessing an equivalent measurable version is presented.
The theorem is  a generalization and correction of  \cite[Theorem 3.8]{mu25}.
\end{abstract}


\date{\today}

\maketitle

\section{Preliminaries}
Throughout this paper $(X,\mfA,P)$ and $(Y,\mfB,Q)$ are fixed probability spaces.  The notation is as presented in the introduction.  If $(Z,\mfZ,T)$ is a measure space, then $(Z,\wh\mfZ,\wh{T})$ denotes its completion and $\mfZ_0:=\{A\in\mfZ:T(A)=0\}$. $T^*$ and $T_*$ are respectively the outer and inner measures induced by $T$.  $\mcL^{\infty}(T)$ is the collection of bounded $T$-measurable real functions.
\begin{deff} \rm Let $R\in{P}\cd{Q}$ be a probability measure . Assume that for every $y\in{Y}$ there exists a $\sigma$-algebra $\mfA_y\subset\mfA$ and a probability $S_y$ on $\mfA_y$ with the following properties:
\begin{description}
\item [(Dis1)]\label{Dis1}
for each $A\in\mfA$ there exists  $N\in\mfB_0$ such that $A\in\mfA_y$ for all $y\notin{N}$ and the function $Y\smi{N}\ni y\longrightarrow S_y(A)$ is $\mfB$-measurable on $Y\smi{N}$;
\item[(Dis2)]
If $A\in\mfA$ and $B\in\mfB$, then
$$
\int_BS_y(A)\,dQ(y)={R}(A\times{B})\,.
$$
\end{description}
The family $\S:=\{(\mfA_y,S_y):y\in{Y}\}$ is called a  $Q$-disintegration of ${R}$ (or a disintegration of $R$ with respect to $Q$).  If $\mfA_y=\mfA$ for all $y\in{Y}$, then $\{(\mfA,S_y):y\in{Y}\}$ is called a  regular conditional probability on $\mfA$ with respect to $Q$.

We say that $P$ is approximated by a family $\mfD\subset\mfA$ (or $P$ is inner regular with respect to $\mfD$)  if for each $A\in\mfA$ and $\ve>0$ there exists $D\in\mfD$ such that $D\subset{A}$ and $P(A\smi{D})<\ve$. $\mfD$ is compact, if $\bigcap_nD_n=\emp$ yields $\bigcap_{n=1}^mD_n=\emp$ for a certain $m\in\N$ and arbitrary $D_n\in\mfD$. $P$ is compact (or countably compact) if it is approximated by a compact family (see \cite{mar}).  \hfill$\Box$
\end{deff}
Disintegration not always exists. The following result of Pachl describes the relation between compactness and existence of disintegration. It is also a partial justification of investigating of behaviour of liftings in product spaces. It proves that the class of spaces admitting disintegrations is very large, in particular it contains all compact measures.
\begin{prop}\label{L2}\cite[Theorem 2.2, Theorem 3.5]{Pachl1}
Let $(X,\mfA,P)$ be a probability space. $P$ is compact if and only if for every complete probability space $(Y,\mfB,Q)$ and every $R\in{P}\cd{Q}$  there exists a $Q$-disintegration of $R$. If $P$ is compact, then  the disintegration $\S=\{(\mfA_y,S_y):y\in{Y}\}$ can be chosen in such a way that $P$ and all measures $S_y$ are approximated also  by the same compact class $\mcK\subset\bigcap_{y\in{Y}}\mfA_y$.
\end{prop}
Let
\begin{align*}
\mcN:=&\{E\subset{X\times{Y}}\colon \forall\;y\in{Y}\;E^y\in\mfA_y\;\exists\;N\in\mfB_0\;\forall\;y\notin{N}\;S_y(E^y)=0\}\\
\wh\mcN:=&\{E\subset{X\times{Y}}\colon \forall\;y\in{Y}\;E^y\in\wh{\mfA_y}\;\exists\;N\in\mfB_0\;\forall\;y\notin{N}\;\wh{S_y}(E^y)=0\}
\end{align*}
Elements of $\mcN$  and $\wh\mcN$ form  $\sigma$-ideals in $\mcP(X\times{Y})$. If $\mfC\subset\mfA$ is a $\sigma$-algebra, then
$\mfC\dd_{\S}\mfB:=\{W\triangle{N}\colon W\in\mfC\otimes\mfB\;\;\&\; N\in\mcN\}$ and $\mfC\dd_{\wh\S}\mfB:=\{W\triangle{N}\colon W\in\mfC\otimes\mfB\;\;\&\; N\in\wh\mcN\}$ are in general larger than $\mfA\otimes\mfB$. If $V=W\tr{N}$ ($W\in\mfC\otimes\mfB$ and $N\in\wh{N}$), then the formula $R_{\dd}(W\triangle{N}):=R(W)$ defines an extension of $R|_{\mfC\otimes\mfB}$.  Obviously, $\mfC\otimes\mfB$ is $R_{\dd}$-dense in $\mfC\dd_{\wh\S}\mfB$. One can  check that $R_{\dd}(V)=\int_YS_y(V^y)\,dQ(y)$ for $V\in\mfC\dd_{\wh\S}\mfB$.
$S_y$-measurability ($\wh{S_y}$-measurability) of  $E^y$ for every element of $E\in\mfC\dd_{\S}\mfB$ ($E\in\mfC\dd_{\wh\S}\mfB$) is the property needed for the main result of the paper.

Because of the condition \eqref{Dis1} one cannot define in the same way the algebra $\mfA\dd_{\S}\mfB$ or $\mfA\dd_{\wh\S}\mfB$. In order to omit the difficulty, we have to assume something more about $\mfC$ and to define a larger algebra. Our general assumption, being a consequence of Proposition \ref{L2},  is the following:
\begin{description}
\item{[C1]}
$\mfC\subset\bigcap_{y\in{Y}}\mfA_y$;
 \item{[C2]}
$\mfC$ is $P$-dense in $\mfA$ and $S_y$-dense in each $\mfA_y$.
\end{description}
Then, let
\[
 \mcM:=\{E\subset{X\times{Y}}\colon \exists\;N\in\mfB_0\;\forall\;y\notin{N}\;\wh{S_y}(E^y)=0\}.
\]
If $\mfA_y=\mfA$ for every $y\in{Y}$, then $\mcM$ is called a family of  $R$-left nil sets  (compare with \cite[3.2.2, Satz 1, Satz 2]{Haupt}). One can also compare with  \cite{BM}, where a similar family of sets (called nil sets) in case of the direct product of measures are defined in a different way. We will keep that terminology also in this more general situation. The measure $R_{\dd}$ can be uniquely extended to a measure on the $\sigma$-algebra
 \[
 \mfA\dd\mfB:=\{V=W\tr{N}\colon W\in\mfA\otimes\mfB\;\&\;N\in\mcM\}.
 \]
 One can check that the extension is the completion of $R_{\dd}$ and it is denoted therefore by $\wh{R_{\dd}}$.   It is defined  by the formula  $\wh{R_{\dd}}(V)=\int_YS_y(V^y)\,dQ(y)$ or, equivalently, by setting $\wh{R_{\dd}}(V)=R(W)$. It is also clear that $\mfA\wh\otimes\mfB$ - the completion of $\mfA\otimes\mfB$ with respect to $R$ - is contained in $\mfA\dd\mfB$ and $\wh{R_{\dd}}$ is an extension of $\wh{R}$.

Clearly $\mfC\dd_{\wh\S}\mfB\subset\mfA\dd\mfB$ and  if $E\in\mfA\dd\mfB$, then there exists $F\in\mfC\dd_{\S}\mfB$ with $\wh{R_{\dd}}(E\tr{F})=0$.

The following fact about skew products is known in case of regular conditional probabilities. Its proof is standard.
\mg{p1}\begin{prop}\label{p1}
Let $(X,\mfA,P)$ and $(Y,\mfB,Q)$ be probability spaces, $R\in{P\cd{Q}}$ and  $\{(\mfA_y,S_y):y\in{Y}\}$ be  a  $Q$-disintegration of ${R}$.
If $f$ is a bounded real valued $\mfA\dd\mfB$-measurable function on $X \times Y$. Then
\begin{itemize}
\item[(a)]
the function $f^y$ is $\wh\mfA_y$-measurable for $Q-a.e. \; y \in {Y}$;
 \item[(b)]
 if $f=0$ $\wh{R_{\dd}}$-a.e., then $f^y=0$ $\wh{S_y}$-a.e. for $Q-a.e. \; y \in {Y}$ ;
\item[(c)]
the function $y\longmapsto \int_Xf^y(x)d\wh{S_y}(x)$ is $\wh\mfB$-measurable;
\item[(d)]
the equality
\mg{e21}\begin{equation}\label{e21}
\int_{X \times Y}f(x,y)d\wh{R_{\dd}}(x,y) =
 \int_Y\int_Xf^y(x)d\wh{S_y}(x)d\wh Q(y)
\end{equation}
holds true.
\end{itemize}
In particular, if $E\in\mfA\dd\mfB$,
then $\wh{R_{\dd}}(E)=\int_Y\wh{S_y}(E^y)\,d\wh{Q}(y)$.
\end{prop}
The following theorem is the key result in our examination of liftings in products.
\mg{T1}\begin{thm} \cite[Theorem 2.1]{mu25}\label{T1}
Assume that $(X,\mfA,P)$  and $(Y,\mfB,Q)$ are probability spaces, $R\in{P}\cd{Q}$ and $\{(\mfA_y,S_y):y\in{Y}\}$ is a $Q$-disintegration of $R$. Let $\mfC$ be a sub-$\sigma$-algebra of $\bigcap_{y\in{Y}}\mfA_y$  and $\mfM:=\mfC\otimes\mfB$.  If  $f \in {\mathcal L}^{\infty}(\wh{R})$ is arbitrary, then each version ${\mathbb E}_{\mfM}(f)$ of the conditional expectation of $f$ with respect to $\mfM$ is such that  $[{\mathbb E}_{\mfM}(f)]^y$ for  $Q$-almost all $y\in{Y}$  is a version of the conditional expectation ${\mathbb E}^y_{\mfC}(f^y)$ of $f^y$ with respect to $\mfC$ and $S_y$. Moreover, the function $(x,y)\longrightarrow {\mathbb E}^y_{\mfC}(f^y)(x)$ is then measurable with respect to $\mfM$ on a set $X\times(Y\smi{N_f})$, where $N_f\in\mfB_0$.
\end{thm}
\begin{rem}\rm
The proof of  \cite[Theorem 2.1]{mu25} remains valid also for $f\in\mcL^{\infty}(\wh{R_{\dd}})$, but I will apply the original version.
\end{rem}
\section{Liftings in product spaces with disintegration}
The notions of lower densities and liftings are used in a little bit more general context than in \cite{it}. More precisely, if $(Z,\mfZ,T)$ is  a probability space and  $\mfD$ is a sub-$\sigma$-algebra of $\mfZ$, then $\delta:\mfD\to\mfZ$ is a lower density if $\delta(\emp)=\emp$, $\delta(Z)=Z$, $\delta(A\cap{B})=\delta(A)\cap\delta(B)$, $T(\delta(A)\tr{A})=0$ and $\delta(A)=\delta(B)$ whenever $A\equiv{B}$. The collection of such lower densities is denoted by $\vt(T|_{\mfD};T|_{\mfZ})$. If $\mfD=\mfZ$, then we write simply $\vt(T)$.
 Since all densities considered in this paper are lower densities I will use the word ``density'' instead of ``lower density''.

If $f\in\mcL^{\infty}(T)$ and  $\mfD\subset\mfZ$ is arbitrary, then ${\mathbb E}_{\mfD}(f)$ denotes the conditional expectation of $f$ with respect to $\mfD$.  If $f\in\mcL^{\infty}(S_y)$, then the conditional expectation of $f$ will be denoted by ${\mathbb E}_{\mfD}^y(f)$.

We begin with a few facts, known for ordinary densities.
\mg{L1}\begin{lem}\label{L1}
Let $(Z,\mfZ,T)$ be a probability space and   $\delta\in\vt(\wh{T})$.  If $\{C_i:i\in{I}\}\subset\mfZ$ is such that $C_i\subset\delta(C_i)$ for every $i\in{I}$,  then $\bigcup_{i\in{I}}C_i\in\wh\mfZ$ and $\bigcup_{i\in{I}}C_i\subset\delta\bigl(\bigcup_{i\in{I}}C_i\bigr)$.
\end{lem}
\begin{lem}\label{L3}\mg{L3}
Let $(Z,\mfZ,T)$ be a probability space and let  $\mfC$ and $\mfD$ be sub-$\sigma$-algebras of $\mfZ$ such that $\mfC\subset\mfD$ and
$\mfD_0\subset\mfC$. Moreover, let $\delta\in\vartheta(T|_\mfC;T|_\mfZ)$ be arbitrary
and let $M\in\mfD\setminus\mfC$.
If $M_1\supset{M}$ and $M_2\supset{M^c}$ are $\mfC$-envelopes of $M$ and
$M^c$ respectively, then the formula
\begin{align*}
 &\widetilde\delta\Bigl[(G \cap M) \cup (H \cap M^c)\Bigr] :=  \\
   &  \Bigl[M \cap \delta\Bigl((G \cap M_1) \cup (H \cap
M^c_1)\Bigr)\Bigr] \cup \Bigl[M^c \cap\delta\Bigl((H \cap M_2) \cup (G
       \cap M^c_2)\Bigr)\Bigr]
\end{align*}
defines a density $\widetilde\delta\in\vartheta(T|_{\sigma(\mfC\cup\{M\}};T|_\mfZ))$
that is an extension of $\delta$.
\end{lem}

\mg{L4}\begin{lem}\label{L4}
Let $(Z,\mfZ,T)$ be a probability space and let  $\mfC$  be a sub-$\sigma$-algebras of $\mfZ$, that is $T$-dense in $\mfZ$. If $\rho:\mfC\to\mfZ$ is a density (lifting), then the formula $\wt\rho(E)=\rho(F)$, where $T(E\tr{F})=0$, defines  a density (lifting) $\wt\rho:\mfZ\to\mfZ$.
\end{lem}
\begin{deff} \rm{ Let $(Z,\mfZ,T)$ be a probability space. A density
$\tau\in\vartheta(T)$ is called an {\em admissible density} (see \cite{mms2}) if it can be
constructed
with the help of the transfinite induction in the way described below.\\
Let $\mfZ = \sigma\{M_{\alpha}:\alpha < \kappa\}$ be
numbered by the ordinals less than $\kappa$, where $\kappa$ is the first ordinal with this property, $M_0=\emp$ and,   for each $1 \leq \gamma < \kappa$ let  $\mfZ_{\gamma}=\sigma(\{M_{\alpha} :\alpha< \gamma \} \cup \mfZ_0)$. Moreover, for each limit ordinal $\gamma<\kappa$, let $\{\gamma_n:n\in\N\}$ be a fixed increasing sequence of ordinals cofinal with $\gamma$.\\
\textbf{(1)}\quad
$\tau_1\in {\vartheta(T|\mfZ_1)}$ is the only existing density on
$(Z,\mfZ_1)$, i.e. $\tau_1(B)=\emp$ if $B\in\mfZ_0$ and, $\tau_1(B)=Z$, if $B\notin\mfZ_0$.\\
\textbf{(2)}\quad
If $\gamma<\kappa$ is a limit ordinal of uncountable cofinality,
then $\mfZ_{\gamma} = \bigcup_{\alpha < \gamma} \mfZ_{\alpha}$ and
we define $\tau_{\gamma}\in\vartheta(T|{\mfZ_{\gamma}})$ by setting
$$
\tau_{\gamma}(B) :\,=
\tau_{\alpha}(B)\quad\mbox{if}\quad B \in {\mfZ_{\alpha}}\quad\mbox{and}
\quad\alpha<\gamma\,.
$$
\textbf{(3)}\quad
Let now $\gamma$ be a limit odinal of countable cofinality.
For simplicity put $\tau_n :\,= \tau_{\gamma_n}$ and
$\mfZ_n :\,= \mfZ_{\gamma_n}$ for all $n \in \N$. Then $\mfZ_
{\gamma} = \sigma(\cup_{n \in \N}\mfZ_n)$ and we can define
$\tau_{\gamma}$ by setting
\begin{equation}\label{e20}
\tau_{\gamma}(B) := \bigcap_{k \in \N}\bigcup_{n \in \N}\bigcap_
{m \geq n}\tau_m(\{{\mathbb E}_{\mfZ_m}(\chi_{{}_B}) > 1-2^{-k} \})\quad\mbox{for}\quad
B \in {\mfZ_{\gamma}}\,.
\end{equation}
It follows from \cite[Lemma 1]{gvw} that $\tau_{\gamma} \in {\vartheta(T|\mfZ_{\gamma})}$ and $\tau_{\gamma} |\mfZ_n = \tau_n$
for each $n\in\N\,.$\\
\textbf{(4)}\quad
Let now $\gamma=\beta+1\,.$
It is well known, that
$$
\mfZ_{\gamma} = \{(G \cap M_{\beta}) \cup (H \cap M_{\beta}^c) : G, \, H \in \mfZ_{\beta}
\}\,.
$$
Let $F_1\supseteq M_{\beta}$ and $F_2\supseteq M_{\beta}^c$ be $\mfZ_{\beta}-$envelopes of
$M_{\beta}$ and $M_{\beta}^c$ respectively. We define a new density by
\begin{align*}
 &\tau_{\gamma}\Bigl((G \cap M_{\beta}) \cup (H \cap M_{\beta}^c)\Bigr) :=  \\
   &  \Bigl(M_{\beta} \cap \tau_{\beta}\Bigl((G \cap F_1) \cup (H \cap
F^c_1)\Bigr)\Bigr) \cup \Bigl(M_{\beta}^c \cap\tau_{\beta}\Bigl((H \cap F_2) \cup (G\cap F^c_2)\Bigr)\Bigr)
\end{align*}
for $G, H \in \mfZ_{\beta}\,.$
By Lemma \ref{L3}  $\tau_{\gamma} \in {\vartheta(T|\mfZ_{\gamma})}$ and
$\tau_{\gamma} |\mfZ_{\beta} = \tau_{\beta}\,.$\\
\textbf{(5)}\quad
We define $\tau\in\vartheta(T)$ by setting $\tau(E)=\tau_{\kappa}(E)$ for $E\in\mfZ$.} 
\end{deff}
\begin{rem}\rm
We might have assumed in the above definition that $M_{\gamma}\notin\mfZ_{\gamma}$ for all $\gamma$ but in Theorem \ref{T2} we need a definition free of such an assumption.
\end{rem}
\begin{deff} \rm Let $\tau\in\vartheta(\mu)$ be an admissible density. If $\rho\in\vL(\mu)$ is such that $\tau(A)\subset\rho(A)$ for every $A\in\mfZ$, then $\rho$ is called lifting  admissibly generated by $\tau$.
\end{deff}
\begin{deff} \rm {Let $(X,\mfA,P)$ be a probability space and $\mfC$
 be a sub-$\sigma$-algebra of $\mfA$. Moreover let  $\S:=\{(\mfA_y,S_y)\colon y\in{Y}\}$ be a disintegration of $R$ with respect to $Q$. Assume that $\mfC:=\sigma(\{M_{\alpha}\colon\alpha<\kappa\})\subset\bigcap_{y\in{Y}}\mfA_y$ is a $\sigma$-algebra and $\mfC_{\gamma}:=\sigma(\{M_{\alpha}:\alpha<\gamma\})$ for $\gamma<\kappa$. We assume that $M_0=\emp$ and $\kappa$ is the smallest ordinal such that $M_{\beta}\notin\mfC_{\alpha}$ is also fulfilled, whenever $\alpha<\beta$.
Moreover, let
\[
\mfC_{y0}:=\{C\in\mfC: S_y(C)=0\}\,,\quad\mbox{ and }
\quad \mfC_{y\gamma}:=\sigma(\mfC_{\gamma}\cup\mfC_{y0})\quad\mbox{for every }y\in{Y}.
\]
For each $y\in{Y}$ we are given an admissible density $\tau_y$ on $\mfC$, that is constructed along the transfinite sequence $\{M_{\alpha}\colon \alpha<\kappa\}$ with fixed  cofinal sequences $\{\gamma_n:n\in\N\}$ and for the transfinite sequence of $\sigma$-algebras $\mfC_{y\beta}$. It is important to notice that the construction of $\tau_y$ yields the inclusion $\tau_y(\mfC_{y\beta})=\tau_{y\beta}(\mfC_{y\beta})\subset \mfC_{y\beta}$ for every $y\in{Y}$.  Therefore in the proofs we will sometimes write $\tau_y$ in place of $\tau_{y\beta}$. We call such a collection $\{\tau_y:y\in{Y}\}$ an equi-admissible family of densities determined by $\{M_{\alpha}\colon \alpha<\kappa\}$.  Due to the density assumption each $\tau_y$ is uniquely extendable to $\wh{\mfA_y}$.}
\end{deff}
\mg{T2}\begin{prop}\label{T2}
Let $(X,\mfA,P)$  and $(Y,\mfB,Q)$ be probability spaces and $R\in{P}\cd{Q}$. Let $\S=\{(\mfA_y,S_y):y\in{Y}\}$ be a disintegration of $R$  with respect to $Q$. Assume that $\mfC\subset\bigcap_{y\in{Y}}\mfA_y$ is such a $\sigma$-algebra that $P$ and each $S_y$ is inner regular with respect to $\mfC$. Assume moreover that $\{\tau_y:y\in{Y}\}$ is an equi-admissible family of densities determined by $\{M_{\alpha}\colon \alpha<\kappa\}$ generating $\mfC$.
Then there
exists a density  $\varphi: \mfA\dd\mfB\to\mfC\dd_{\S}\mfB$  such that for
each $F\in{\mfA}\dd\mfB$
\begin{equation}\label{e18}
[\varphi(F)]^y = \tau_y([\varphi(F)]^y)\qquad\mbox{for  all }\;
y \in Y\,.
\end{equation}
\end{prop}
\begin{proof}  First we will construct a density   $\varphi: \mfC\otimes\mfB\to\mfC\dd_{\S}\mfB$.
To do it, we will use the transfinite induction.

To simplify further notation I denote each $\sigma$-algebra $\mfC_{\alpha}\dd_{\S}\mfB:=\sigma({\mfC}_{\alpha}\otimes \mfB\cup\mcN)$ by $\mfN_{\alpha}$ and the $\sigma$-algebra $\sigma((\mfC_{\alpha}\otimes\mfB)\cup(\mfC\otimes\mfB)_0)$ by $\mfM_{\alpha}$.  Let $\delta\in\vt(Q)$ be a fixed density.
At the first step we  define $\varphi_1\in\vartheta\Bigl(R_{\dd}|_{\mfM_1}\Bigr)$ by
$\varphi_1(E)=\emp$ if $R_{\dd}(E)=0$ and $\varphi_1(E)=X\times\delta(B)$, if $E\equiv_{R_{\dd}}X\times{B}$ with $B\in\mfB$.

Let us fix now $\gamma<\kappa$ and assume  that    for each $\alpha<\gamma$  there exists a density
$\varphi_{\alpha}\in\vt(R_*|{{}_{\mfM_{\alpha}}};R_*|_{{}_{\mfN_{\alpha}}})$
such that for each $F \in \mfM_{\alpha}$
\begin{equation}\label{e4}
[\varphi_{\alpha}(F)]^y\in\mfC_{y\alpha}\quad\mbox{and}\quad[\varphi_{\alpha}(F)]^y=\tau_y([\varphi_{\alpha}(F)]^y)\qquad\mbox{for every}\; y \in Y
\end{equation}
and
\begin{equation}\label{e3}
\varphi_{\beta}(F)=\varphi_{\alpha}(F)\quad \mbox{if }\alpha<\beta<\gamma\;\mbox{and\;}F\in\mfM_{\alpha}.
\end{equation}
We will split now the proof into three parts:\\
\textbf{(A)} $\gamma$ is a limit ordinal of countable cofinality. \\
Let $(\gamma_m)_m$ be the increasing sequence of ordinals cofinal with $\gamma$.
Applying \eqref{e20} we define $\ov\varphi_{\gamma}\in\vt(R_*|{{}_{\mfM_{\gamma}}};R_*|_{{}_{\mfN_{\gamma}}})$ extending all $\varphi_{\gamma_m}$, setting for each $F\in\mfM_{\gamma}$
\begin{equation*}
\ov\varphi_{\gamma}(F):= \bigcap_{k \in \N}\bigcup_{n\in\N}
\bigcap_{m>n}\varphi_{\gamma_m}\left(\left\{{\mathbb E}_{\mfM_{\gamma_m}}(\chi_{{}_F})
> 1 - 2^{-k}\right\}\right).
\end{equation*}
According to \eqref{e20} we have also for each $A \in {\mfC_{\gamma}}$
\begin{equation*}
\tau_y(A):= \bigcap_{k \in \N}\bigcup_{n\in\N}
\bigcap_{m>n}\tau_y\left(\left\{{\mathbb E}^y_{\mfC_{\gamma_m}}(\chi_{{}_A}) > 1 -
2^{-k}\right\}\right)
\end{equation*}
 Then for $F \in \mfM_{\gamma}$ let $N_{m,F}\in\mfB_0$ be such that (see
Theorem \ref{T1}) for each $y\notin{N_{m,F}}$ .
$$
[{\mathbb E}_{\mfC_{\gamma_m} \otimes {\mathfrak
B}}(\chi_{{}_F})]^y={\mathbb E}^y_{\mfC_{\gamma_m}}(\chi_{{}_{F^y}})\;\;
S_y|\mfC_{\gamma_m}-a.e.
$$
and (by the inductive assumption)
$$
[\varphi_{\gamma_m}(\{{\mathbb E}_{\mfC_{\gamma_m} \otimes \mfB}
       (\chi_{{}_F}) >1 - 2^{-k}\})]^y=\tau_y([\{{\mathbb E}_{\mfC_{\gamma_m} \otimes \mfB}
       (\chi_{{}_F}) >1 - 2^{-k}\}]^y).
$$
If $ N_F:= \bigcup_{m \in\N} N_{m,F}$, then $N_F \in
\mfB_0$ and it follows from Theorem \ref{T1}, that
for all $y \in Y \setminus N_F$
\begin{align*}
    [\ov\varphi_{\gamma}(F)]^y & =  \bigcap_{k \in \N}\bigcup_{n \in \N}
      \bigcap_{m \geq n}\left[\varphi_{\gamma_m}\left(\left\{{\mathbb E}_{\mfC_{\gamma_m} \otimes \mfB}
       (\chi_{{}_F}) >1 - 2^{-k}\right\}\right)\right]^y \\
        & =  \bigcap_{k \in \N}\bigcup_{n \in \N}\bigcap_{m
              \geq n}\tau_y\left(\left[\left\{{\mathbb E}_{\mfC_{\gamma_m} \otimes \mfB}
              (\chi_{{}_F}) > 1 - 2^{-k}\right\}\right]^y\right) \\
        & =  \bigcap_{k \in \N}\bigcup_{n \in \N}\bigcap_{m
              \geq n}\tau_y\left(\left\{[{\mathbb E}_{\mfC_{\gamma_m} \otimes \mfB}
              (\chi_{{}_F})]^y > 1 - 2^{-k}\right\}\right) \\
        & =  \bigcap_{k \in \N}\bigcup_{n \in\N}\bigcap_{m \geq n}
              \tau_y\left(\left\{{\mathbb E}^y_{\mfC_{\gamma_m}}(\chi_{{}_{F^y}}) > 1 - 2^{-k}\right\}\right) =
              \tau_y(F^y).
\end{align*}
It follows that $\tau_y([\ov\varphi_{\gamma}(F)]^y)=[\ov\varphi_{\gamma}(F)]^y$,
for all $y \notin M_F\in\mfB_0$, where $N_F\subset{M_F}$.  Thus, $D:=\{(x,y):\tau_y([\ov\varphi_{\gamma}(F)]^y)\neq[\ov\varphi_{\gamma}(F)]^y\}\in\mcN$, because if $y\in{M_F}$, then $[\ov\varphi_{\gamma}(F)]^y\in\mfC_{y\gamma}$. It follows that we may define $\varphi_{\gamma}$ by setting $[\varphi_{\gamma}(F)]^y :=\tau_y([\ov\varphi_{\gamma}(F)]^y)$ for every $y\in{Y}$.  Obviously $\varphi_{\gamma}(F)\in\mfN_{\gamma}$. One can also see that if $F\in\mfM_{\alpha}$ for some $\alpha<\gamma$, then there is $\gamma_m>\alpha$ such that $F\in\mfM_{\gamma_m}$, what yields $\ov\varphi_{\gamma}(F)=\varphi_{\gamma_m}(F)=\varphi_{\alpha}(F)$. Hence $[\varphi_{\gamma}(F)]^y=
\tau_{y\gamma}\bigl([\ov\varphi_{\gamma}(F)]^y\bigr)=\tau_{y\gamma}\bigl([\varphi_{\alpha}(F)]^y\bigr)=[\varphi_{\alpha}(F)]^y$.\\
\textbf{(B)} $\gamma=\beta+1$. \\
 Let $W_1\supset{M_{\beta}\times{Y}}$ and $W_2\supset{M_{\beta}^c\times{Y}}$ be respectively $\mfM_{\beta}$-envelopes of $M_{\beta}\times{Y}$ and $M_{\beta}^c\times{Y}$ with respect to $R_{\dd}|_{\mfM_{\beta}}$.
We define an extension $\ov\varphi_{\gamma}:\mfM_{\gamma}\to \mfN_{\gamma}$ of $\varphi_{\beta}$  by the formula (see Lemma \ref{L3})
\begin{align*}
 &\ov\varphi_{\gamma}\Bigl[\bigl(G \cap (M_{\beta}\times{Y})\bigr) \cup \bigl(H \cap (M_{\beta}^c\times{Y})\bigr)\Bigr] :=  \\
   &  \Bigl[(M_{\beta}\times{Y}) \cap \varphi_{\beta}\Bigl((G \cap W_1) \cup (H \cap
W^c_1)\Bigr)\Bigr] \cup \Bigl[(M_{\beta}^c\times{Y}) \cap\varphi_{\beta}\Bigl((H \cap W_2) \cup (G
       \cap W^c_2)\Bigr)\Bigr]
\end{align*}
if $G,H\in\mfM_{\beta}$.

It follows from  the basic properties of densities that  for each $F\in\mfM_{\gamma}$ the equality $R_{\dd}(\ov\varphi_{\gamma}(F)\tr{F})=0$ is valid. In virtue of Proposition \ref{p1} there exists a set $N\in\mfB_0$ such that for every $y\notin{N}$ we have $\tau_{y\gamma}([\ov\varphi_{\gamma}(F)]^y)=\tau_{y\gamma}(F^y)$. Consequently,   if the set $N$ is  chosen  for $W:=\bigl(G \cap (M_{\beta}\times{Y})\bigr) \cup \bigl(H \cap (M_{\beta}^c\times{Y})\bigr)$,  $(G \cap W_1) \cup (H \cap{W^c_1})$ and $(H \cap W_2) \cup (G
       \cap {W^c_2})$, then we have for $y\notin{N}$
\[
 \mfC_{y\gamma}\ni[\ov\varphi_{\gamma}(W)]^y \equiv_{{}_{S_y}}\tau_{y\gamma}\bigl([\ov\varphi_{\gamma}(W)]^y\bigr)\in\mfC_{y\gamma}.
\]
If $y\in{N}$, then $[\ov\varphi_{\gamma}(W)]^y\in\mfC_{y\gamma}$ as well. As a result, if $D^y:=[\ov\varphi_{\gamma}(W)]^y \tr\tau_{y\gamma}\bigl([\ov\varphi_{\gamma}(W)]^y\bigr)$, then $D\in\mcN$ and so we may  define
$\varphi_{\gamma}:\mfM_{\gamma}\to\mfN_{\gamma}$ setting for every $y\in{Y}$
 $[\varphi_{\gamma}(W)]^y:=\tau_{y\gamma}([\ov\varphi_{\gamma}(W)]^y)$. One can easily check that $\ov\varphi_{\gamma}|_{{\mfM}_{\beta}}=\varphi_{\beta}$ and so if $E\in{\mfM}_{\beta}$, then $[\varphi_{\gamma}(E)]^y=\tau_{y\gamma}([\ov\varphi_{\gamma}(E)]^y)=\tau_{y\gamma}([\varphi_{\beta}(E)]^y)=[\varphi_{\beta}(E)]^y$. It follows that  $\varphi_{\gamma}|_{{\mfM}_{\beta}}=\varphi_{\beta}$.\\
\textbf{(C)} Assume now that $\gamma$ is of uncountable cofinality. \\
In such a case $\mfC_{\gamma}=\bigcup_{\alpha<\gamma}\mfC_{\alpha}$ and for each $F\in\mfC_{\gamma}$ there exists $\alpha<\gamma$ with $F\in\mfC_{\alpha}$. We simply set $\varphi_{\gamma}(F)=\varphi_{\alpha}(F)$.\\
\textbf{(D)}  When $\gamma=\kappa$, then we get a density $\varphi_{\kappa}\in\vartheta(R_*|_{\mfM_{\kappa}};R_*|_{\mfN_{\kappa}})$ and the inductive part is finished.\\
\textbf{(E)} Now we are going to extend $\varphi_{\kappa}$ to a density $\varphi:\mfA\dd\mfB\to \mfC\dd_{\S}\mfB=\mfN_{\kappa}$.  As it has been mentioned earlier $\mfC\dd_{\S}\mfB$ is $\wh{R_{\dd}}$-dense in  $\mfA\dd\mfB$ and $\mfC\otimes\mfB$ is  $\wh{R_{\dd}}$-dense in $\mfC\dd_{\S}\mfB$.  If $W\in\mfA\dd\mfB$ and $V\in\mfC\otimes\mfB$ us such that $\wh{R_{\dd}}(W\tr{V})=0$, then we set $\varphi(W):=\varphi_{\kappa}(V)$.
$\varphi:\mfA\dd\mfB\to \mfC\dd_{\S}\mfB$ and the equality (\ref{e18}) holds true.
\end{proof}
\begin{prop}\label{p3}
Under the assumptions of Theorem \ref{T2}
there exists a density $\psi:\mfA\dd\mfB\longrightarrow\mfC\dd_{\wh\S}\mfB$ satisfying for all $E\in\mfA\dd\mfB$ the following conditions:
\begin{itemize}
\item[(1)]
$ \wh{S_y}\Bigl([\psi(E)]^y\cup[\psi(E^c)]^y\Bigr)=1 \qquad\mbox{for
all}\quad y\in Y\,;$
\item[(2)]
$[\psi(E)]^y= \tau_y\left([\psi(E)]^y\right) \qquad\mbox{for all}
\quad  y\in Y\,;$
\end{itemize}
\end{prop}
\begin{proof}  For each $y\in{Y}$ the measure $S_y$ is inner regular with respect to $\mfC$ and consequently the density $\tau_y\in\vartheta(S_y|_{\mfC})$ can be uniquely extended to $\wh{S_y}|_{\mfC_y}$. We denote the extension also by $\tau_y$. The density $\varphi: \mfA\dd\mfB\to \mfC\dd_{\S}\mfB$ satisfying the thesis of Proposition \ref{T2} can be treated as an element of $\vt(\wh{R_{\dd}})$.  Let
\begin{equation*}
\Psi:\,=\left\{\wt\varphi\in\vartheta(\wh{R_{\dd}})\;\left|\
\begin{aligned}
&
\forall E\in\mfA\dd\mfB\;
\varphi(E)\subseteq\wt\varphi(E)\;\&\;\wt\varphi(E)\in\mfC\dd_{\wh\S}\mfB\, \ \\[-1mm]
&\;\&\;\forall y\in{Y}\;[\wt\varphi(E)]^y\subseteq\tau_y([\wt\varphi(E)]^y)\,
\end{aligned}
\right.\right\}\,.
\end{equation*}
We consider $\Psi$ with inclusion as the partial order:
$\wt\varphi_1\leq\wt\varphi_2$ if
$\wt\varphi_1(E)\subseteq\wt\varphi_2(E)$ for each
$E\in\mfA\dd\mfB$.

\textbf{(I)} {\it There exists a maximal element in $\Psi$.}

The only fact we have to prove is showing that each chain
$\{\varphi_{\alpha}\}_{\alpha\in A}\subseteq\Psi$ has a dominating
element in $\Psi$. The obvious candidate is $\wt\varphi$
given for each $E\in\mfA\dd\mfB$ by
\begin{equation}\label{xx}
\wt\varphi(E)=\bigcup_{\alpha\in A}\varphi_{\alpha}(E)\,.
\end{equation}
One can easily check that $\wt\varphi(E)\subseteq[\wt\varphi(E^c)]^c$ and
so if an $\alpha\in A$ is fixed, then
$$
\varphi_{\alpha}(E)\subseteq\wt\varphi(E)\subseteq[\wt\varphi(E^c)]^c
\subseteq[\varphi_{\alpha}(E^c)]^c
$$
for each $E\in\mfA\dd\mfB$. Since $\wh{R_*}$ is complete and $\varphi_{\alpha}\in\vartheta(
\wh{R_*})$, this proves that  $\wt\varphi(E)\in\mfA\dd\mfB$
and $\wt\varphi(E)\stackrel{\wh{R_*}}{=}E$.
Consider now the section properties of $\wt\varphi(E)$. For
fixed $ y\in Y$
$$
[\wt\varphi(E)]^y=\bigcup_{\alpha\in
A}[\varphi_{\alpha}(E)]^y \subseteq\bigcup_{\alpha\in
A}\tau_y([\varphi_{\alpha}(E)]^y)
$$
and so - in virtue of Lemma \ref{L1} - the set
$[\wt\varphi(E)]^y$ is $\wh{S_y}$-measurable. As a result $\wt\varphi(E)\in\mfC\dd_{\wh\S}\mfB$. Setting in the
above inclusions $\wt\varphi(E)$ instead of
$\varphi_{\alpha}(E)$, we see that the inclusion
$[\wt\varphi(E)]^y\subseteq\tau_y([\wt\varphi(E)]^y)$
is satisfied also.

The multiplication rule and other properties are clear and so
$\wt\varphi\in\vartheta(\wh{R_*})$ and  $\varphi:\mfA\dd\mfB\longrightarrow\mfA\dd_{\wh\S}\mfB$. That proves that
$\wt\varphi$ dominates the whole chain.

According to the Zorn-Kuratowski Lemma the set $\Psi$ possesses a
maximal element $\psi:\mfA\dd\mfB\longrightarrow\mfC\dd_{\wh\S}\mfB$.

\textbf{(II)}  {\it For each $ y\in Y$ and $E\in\mfA\dd\mfB$}\; $\wh{S_y}([\psi(E)]^y\cup[\psi(E^c)]^y)=1$.

 Notice first that as $\varphi\in\Psi$ all sections $[\psi(E)]^y$ are $\wh{S_y}$-measurable.
Suppose now that there exist $H\in\mfA\dd\mfB$ and $ y_0\in Y$ such that $
\wh{S}_{y_0}([\psi(H)]^{ y_0}\cup[\psi(H^c)]^{ y_0})<1\,.$ Let
$$
W:\,=\tau_{y_0}\Bigl[\Bigl([\psi(H)]^{ y_0}\cup[\psi(H^c)]^{ y_0}
\Bigr)^c\Bigr]
$$
and
\[
[\wh\psi(E)]^y=\left\{\begin{array}{ll}
[\psi(E)]^y&$ if $ y\neq y_0\\

[\psi(E)]^{y_0}\cup(W\cap[\psi(H\cup E)]^{y_0}) &$if $\; y= y_0
\end{array}\right.
\]
It is clear, that $\psi(E)\subseteq\wh\psi(E)$ for each $E\in
\mfA\dd\mfB\,.$ In particular $\wh\psi(
X\times Y)= X\times Y\,.$ Also the other density properties are
fulfilled.

It follows directly from the definition that $\wh\psi$ and $\psi$
are different densities. In order to get a contradiction with our
hypothesis it is enough to show that $[\wh\psi(E)]^{
y_0}\subseteq\tau_{y_0}\Bigl([\wh\psi(E)]^{ y_0}\Bigr)$, but this
is immediate.
This completes the proof.

\textbf{(III)}  {\it For each $ y\in Y$ and $E\in\mfA\dd \mfB$ \;$[\psi(E)]^y=\tau_y\Bigl([\psi(E)]^y\Bigr)$. }

Set for each $ y\in Y\;\mbox{and}\; E\in\mfA\dd\mfB$
$$
[\widetilde\psi(E)]^y=\tau_y\Bigl([\psi(E)]^y\Bigr)\,.
$$
Clearly $\psi(F)\subseteq\widetilde\psi(F)$ for each $F$. Moreover
the equality $\psi(E)\cap\psi(E^c)=\emptyset$ yields for each $y$
the relation $\tau_y\Bigl([\psi(E)]^y\Bigr)\cap
\tau_y\Bigl([\psi(E^c)]^y\Bigr)=\emptyset\,.$ As a consequence, we
get $\widetilde\psi(E)\cap\widetilde\psi(E^c)=\emptyset$, and then
$\widetilde\psi(E^c)\subseteq(\widetilde\psi(E))^c\,.$ Hence
$$
\psi(E^c)\subseteq\widetilde\psi(E^c)\subseteq[\widetilde\psi(E)]^c\subseteq
[\psi(E)]^c\,.
$$
Since $\psi$ is a density, we have $R_{\dd}([\psi(E)]^c)=R_{\dd}[\psi(E^c)]$
and so  $\widetilde\psi(E)$ is $\mfA\dd\mfB$-measurable. But at the same time $\wh{S_y}\bigl([\psi(E)]^y\bigr)=\wh{S_y}\bigl([\psi(E^c)]^y\bigr)$. Hence $\wt\psi(E)\in\mfC\dd_{\wh{\S}}\mfB$.  It follows
that $\widetilde\psi\in\Psi$ and so $\psi=\widetilde\psi$, due to
the maximality of $\psi$.
\end{proof}
\mg{T3}\begin{thm}\label{T3}
Let $(X,\mfA,P)$  and $(Y,\mfB,Q)$ be probability spaces. Assume that $R\in{P}\cd{Q}$ has a disintegration ${\mathbb S}=\{(\mfA_y,S_y):y\in{Y}\}$ with respect to $Q$ and $\mfC\subset\bigcap_{y\in{Y}}\mfA_y$ is such a $\sigma$-algebra that $P$ and each $S_y$ are inner regular with respect to $\mfC$. Let $\{\tau_y\in\vartheta(\wh{S_y}):y\in{Y}\}$ be the family from Proposition \ref{p3}.  Then there exist liftings $\sigma_y\in\vL(\wh{S_y})$ that are admissibly generated by $\tau_y$ for all $y\in Y$ and there exists a lifting $\pi:\mfA\dd\mfB\to \mfC\dd_{\wh\S}\mfB$ such that the following condition is satisfied:
\begin{align*}
[\pi(E)]^y&= \sigma_y\Bigl([\pi(E)]^y\Bigr)\;\;  \quad\mbox{for every } E\in \mfA\dd\mfB\;\mbox{and }y\in{Y};\\
[\pi(f)]^y &= \sigma_y\Bigl([\pi(f)]^y\Bigr)\quad\mbox{for every }f \in {\mathcal L}^{\infty}(\wh{R_{\dd}})\;\mbox{and } y\in{Y}.
\end{align*}
\end{thm}
\begin{proof}  We take $\psi $ from Proposition \ref{p3} and liftings $\sigma_y\in\vL(\wh{S_y})$ generated by $\tau_y$. Then we define
$\pi:\mfA\dd\mfB\to \mfC\dd_{\wh\S}\mfB$ by setting for each $E\in \mfA\dd\mfB$ and each $y\in Y$
\begin{equation*}
[\pi(E)]^y= \sigma_y\Bigl([\psi(E)]^y\Bigr)\,.
\end{equation*}
\end{proof}
The next theorem is a direct consequence of Theorem \ref{T3} and it is a strong generalization of \cite[Theorem 3.6]{smm}, where the $\sigma$-algebra $\mfA$ was assumed to be separable in the  Frechet-Nikod\'{y}m pseudometric.
\mg{T4}\begin{thm}\label{T4}
Let $(X,\mfA,P)$  and $(Y,\mfB,Q)$ be probability spaces and  $R\in{P}\cd{Q}$ be arbitrary.  Assume that   $\{(\mfA,S_y):y\in{Y}\}$ is a regular conditional probability on $\mfA$ with respect to $Q$  and  $\{\tau_y\in\vartheta(S_y):y\in{Y}\}$ is an equi-admissible family of densities. Then for each $y\in{Y}$ there exists $\sigma_y\in\vL(\wh{S_y})$ admissibly generated by $\tau_y$ and there exists a lifting  $\pi:\mfA\dd\mfB\to \mfA\dd_{\wh\S}\mfB$  such that the following conditions are satisfied:
\begin{align*}
[\pi(E)]^y&= \sigma_y\Bigl([\pi(E)]^y\Bigr)\;\;  \quad\mbox{for every } E\in \mfA\dd\mfB\;\mbox{and }y\in{Y};\\
[\pi(f)]^y &= \sigma_y\Bigl([\pi(f)]^y\Bigr)\quad\mbox{for every }f \in {\mathcal L}^{\infty}(\wh{R_{\dd}})\;\mbox{and } y\in{Y}.
\end{align*}
\end{thm}
Theorem \ref{T4} is also a partial generalization of the main results from \cite{mms1} and \cite{mms2}. In \cite{mms2} the measure $R$ was a direct product of $P$ and $Q$ whereas in \cite{mms1} the measure $R$ was assumed to be absolutely continuous with respect to the direct  product $P\times{Q}$.
\section{Connections with earlier results.}
\subsection{Absolutely continuous disintegrations}
The topic investigated in this paper was earlier undertaken in \cite{mms2,smm,mms1,mu25}. In \cite{mms2} the measure $R$ was simply the direct product of $P$ and $Q$.  In \cite{smm} the algebra $\mfA$ was assumed to be separable in the Frechet-Nikod\'{y}m metric. In \cite{mms1} the measure $R$ was assumed to be absolutely continuous with respect to $2^{-k}$.  In all the above works the final lifting $\pi$ was an element of $\vL(\wh{R})$ but ${\mathbb S}$ was assumed to be a regular conditional probability. In \cite{mu25} $R$ was arbitrary and ${\mathbb S}$ was a disintegration , but the proof of  the inclusion $\pi\in\vL(\wh{R})$  was erroneus. Theorem \ref{T3} corrects the mistake, proving that in general  $\pi\in\vL(\wh{R_{\dd}})$.

Thus, investigating densities and liftings in the context of disintegration or regular conditional probabilities one immediately meets the following problem:

When the density $\varphi$ in Proposition \ref{T2} may take values in $\mfA\otimes\mfB$ and  the lifting $\pi$ from Theorem \ref{T3} or \ref{T4} may take values in $\mfA\wh\otimes\mfB$\,?

It follows from the proof of Proposition \ref{T2}, that in general such a relation is impossible.

 Below I prove that if $R\ll{P}\otimes{Q}$  in case of the regular conditional probability (or $S_y\ll{P|{_{\mfA_y}}}$ in case of a disintegration), then the density $\varphi$ in Proposition \ref{T2}  and  the lifting $\pi$ from Theorem \ref{T3} or \ref{T4} may take values in $\mfA\wh\otimes\mfB$. That is an alternate proof of the corresponding part of the proof of \cite[Theorem 3.3]{mms1}.

In section (B) of the proof of Proposition \ref{T2}, the density $\ov\varphi_{\gamma}$ was defined by the formula
\begin{align*}
 &\ov\varphi_{\gamma}\Bigl[\bigl(G \cap (M_{\beta}\times{Y})\bigr) \cup \bigl(H \cap (M_{\beta}^c\times{Y})\bigr)\Bigr] :=  \\
   &  \Bigl[(M_{\beta}\times{Y}) \cap \varphi_{\beta}\Bigl((G \cap W_1) \cup (H \cap
W^c_1)\Bigr)\Bigr] \cup \Bigl[(M_{\beta}^c\times{Y}) \cap\varphi_{\beta}\Bigl((H \cap W_2) \cup (G
       \cap W^c_2)\Bigr)\Bigr]\,,
\end{align*}
where $G,H\in\mfM_{\beta}$, $\gamma=\beta+1$ and  $W_1\supset{M_{\beta}\times{Y}}$ and $W_2\supset{M_{\beta}^c\times{Y}}$ are respectively $\mfC_{\beta}\dd_{\S}\mfB$-envelopes of $M_{\beta}\times{Y}$ and $M_{\beta}^c\times{Y}$ with respect to $R_{\dd}$. In particular,
\[
\ov\varphi_{\gamma}(M_{\beta}\times{Y})=
  \Bigl[(M_{\beta}\times{Y}) \cap \varphi_{\beta}(W_1)\Bigr]\cup \Bigl[(M_{\beta}^c\times{Y}) \cap\varphi_{\beta}(  W^c_2)\Bigr]
\]
If $\varphi_{\beta}:\mfM_{\beta}\longrightarrow\mfA\wh\otimes\mfB$, then we have also $\ov\varphi_{\gamma}:\mfM_{\gamma}\longrightarrow\mfA\wh\otimes\mfB$.

We know that there exists $N\in\mfB_0$ such that $\tau_y\bigl(\ov\varphi_{\gamma}(M_{\beta}\times{Y})\bigr)=\tau_y(M_{\beta})$, for every $y\notin{N}$. Assume that $\varphi_{\beta}(\mfA_{\beta}\dd_{\S}\mfB)\subset\mfA\wh\otimes\mfB$. If $\varphi_{\gamma}$ is defined (as in (A)) by $[\varphi_{\gamma}(E)]^y=\tau_y\bigl([\ov\varphi_{\gamma}(E)]^y\bigr)$, then $\varphi_{\gamma}(M_{\beta}\times{Y})\in\mfA\wh\otimes\mfB$ if and only if $\bigcup_y\tau_y(M_{\beta})\times\{y\}\in\mfA\wh\otimes\mfB$. Thus, if there exists an equi-admissible family  $\{\tau_y:y\in{Y}\}$ with the property  $\bigcup_y\tau_y(A)\times{Y}\in\mfA\wh\otimes\mfB$, for every $A\in\mfA$, then we might have  $\varphi_{\gamma}(\mfA_{\gamma}\dd_{\S}\mfB)\subset\mfA\wh\otimes\mfB$, for every $\gamma<\kappa$.

It is so if $R=P\times{Q}$, because then all $\tau_y$'s may coincide with a fixed density on $\mfA$.

In the general case the situation is much more complicated.
Let $V_{1y}\supset{M_{\beta}}$ and $V_{2y}\supset{M_{\beta}^c}$ be respectively $\mfC_{y\beta}$-envelopes of $M_{\beta}$ and $M_{\beta}^c$ with respect to $S_y$. Without loss of generality, we may assume that
$$
V_{1y}\subset{W_1^y}\,,\quad V_{2y}\subset{W_2^y}\quad\mbox{for all }y\in{Y}\,.
$$
The construction of the densities $\tau_y$ guarantees the equality
\mg{e11}\begin{align}\label{e11}
 &\tau_{y\gamma}\Bigl[(A \cap M_{\beta}) \cup (B \cap M_{\beta}^c)\Bigr] :=  \\
   & \Bigl[M_{\beta} \cap \tau_{y\beta}\Bigl((A \cap V_{1y}) \cup (B \cap
V^c_{1y})\Bigr)\Bigr] \cup \Bigl[M_{\beta}^c \cap\tau_{y\beta}\Bigl((B \cap V_{2y}) \cup (A
       \cap V^c_{2y})\Bigr)\Bigr]\notag
\end{align}
if $A,B\in\mfC_{\beta}$. But if  $S_y\ll{P|_{\mfC_y}}$ for every $y\in{Y}$ (If $R\ll{P}\times{Q}$, then such absolutely continuous $\S$ exists), then one can easily check that one may set $V_{1y}=W_1^y$ and $V_{2y}=W_2^y$.

It follows from  the basic properties of densities that for each $F\in\mfM_{\beta}$ there exists a set $N_1\in\mfB_0$ such that for every $y\notin{N_1}$ we have $\tau_y([\varphi_{\beta}(F)]^y)=\tau_y(F^y)$. Similarly, for each $F\in\mfM_{\gamma}$ there exists a set $N_2\in\mfB_0$ such that for every $y\notin{N_2}$ we have $\tau_y([\ov\varphi_{\gamma}(F)]^y)=\tau_y(F^y)$. Consequently, it follows from (\ref{e4}) and \eqref{e11} that  if $N\in\mfB_0$ is chosen for $(G \cap W_1) \cup (H \cap{W^c_1})$, \;$(H \cap W_2) \cup (G
       \cap {W^c_2})$ and $\bigl(G \cap (M_{\beta}\times{Y})\bigr)\cup \bigl(H\cap (M_{\beta}^c\times{Y})\bigr)$, then for $y\notin{N}$
\begin{align*}
 &\biggl[\ov\varphi_{\gamma}\Bigl[\bigl(G \cap (M_{\beta}\times{Y})\bigr) \cup \bigl(H \cap (M_{\beta}^c\times{Y})\bigr)\Bigr]\biggr]^y  \\
  &=\Bigl[M_{\beta} \cap \tau_y\Bigl((G^y \cap W_1^y) \cup (H^y \cap
(W_1^y)^c)\Bigr)\Bigr] \cup \Bigl[M_{\beta}^c \cap\tau_y\Bigl((H^y \cap W_2^y) \cup (G^y
       \cap (W_2^y)^c)\Bigr)\Bigr]\\
       &=  \Bigl[M_{\beta} \cap \tau_y\Bigl((G^y \cap V_{1y}) \cup (H^y \cap
V^c_{1y})\Bigr)\Bigr] \cup \Bigl[M_{\beta}^c \cap\tau_y\Bigl((H^y \cap V_{2y}) \cup (G^y
       \cap V_{2y}^c)\Bigr)\Bigr]\\
       &=\tau_y\Bigl((G^y\cap{M_{\beta}})\cup (H^y\cap{M_{\beta}^c})\Bigr)= \tau_y\biggl[\Bigl[\bigl(G \cap (M_{\beta}\times{Y})\bigr) \cup \bigl(H \cap (M_{\beta}^c\times{Y})\bigr)\Bigr]^y\biggr]\\
         &=\tau_y\biggl(\biggl[\ov\varphi_{\gamma}\Bigl[\bigl(G \cap (M_{\beta}\times{Y})\bigr) \cup \bigl(H \cap (M_{\beta}^c\times{Y})\bigr)\Bigr]\biggr]^y\biggr).
\end{align*}
If $y\in{N}$, then $\biggl[\ov\varphi_{\gamma}\Bigl[\bigl(G \cap (M_{\beta}\times{Y})\bigr) \cup \bigl(H \cap (M_{\beta}^c\times{Y})\bigr)\Bigr]\biggr]^y\in\mfC_{y\gamma}$ and so we may  define
$\varphi_{\gamma}:{\mfM}_{\gamma}\longrightarrow\sigma\bigl(\mfC_{\gamma}\otimes\mfB\cup(\mfA\wh\otimes\mfB)_0\bigr)$ setting for every $y\in{Y}$
 $[\varphi_{\gamma}(W)]^y=\tau_y([\ov\varphi_{\gamma}(W)]^y)$ for every $W\in{\mfM}_{\gamma}$. The family $(\mfA\wh\otimes\mfB)_0$ is necessary, because $X\times{N}\in(\mfA\otimes\mfB)_0$ and  $\varphi_{\gamma}(W)\cap(X\times{N})\neq\emp$ for some $W$.  One can easily check that $\ov\varphi_{\gamma}|_{\mfM_{\beta}}=\ov\varphi_{\beta}$ (defined one step earlier by limit or non-limit procedure) and so $\varphi_{\gamma}|_{{}_{\mfM_{\beta}}}=\varphi_{\beta}$.

If $\gamma$ is a limit ordinal of countable cofinality, then the proof is exactly as in Proposition \ref{T2}. The same concerns the rest of the proof.  \hfill$\Box$

The above considerations can be formulated as
\mg{T5}\begin{thm}\label{T5}
Let $(X,\mfA,P)$  and $(Y,\mfB,Q)$ be probability spaces. Assume that $R\in{P}\cd{Q}$ has a  disintegration ${\mathbb S}=\{(\mfA_y,S_y):y\in{Y}\}$ with respect to $Q$ that is absolutely continuous (i.e. $S_y\ll{P|{\mfA_y}}$ for every $y\in{Y}$) and $\mfC\subset\bigcap_{y\in{Y}}\mfA_y$ is such a $\sigma$-algebra that $P$ and each $S_y$ are inner regular with respect to $\mfC$. Let $\{\tau_y\in\vartheta(\wh{S_y}):y\in{Y}\}$ be the family from Proposition \ref{p3}.  Then there exist liftings $\sigma_y\in\vL(\wh{S_y})$ that are admissibly generated by $\tau_y$ for all $y\in Y$ and there exists a lifting $\pi:\mfA\dd\mfB\to \mfC\wh\otimes\mfB\subset\mfA\wh\otimes\mfB$ such that the following condition is satisfied:
\begin{equation*}
[\pi(f)]^y = \sigma_y\Bigl([\pi(f)]^y\Bigr)\quad \mbox{for every }f \in {\mathcal L}^{\infty}(\wh{R_{\dd}}) \;\mbox{and every }y\in Y.
\end{equation*}
\end{thm}
\subsection{Measurable modifications of stochastic processes}
Let $(X,\mfA,P),\,(Y,\mfB,Q),\;R,\;R_{\dd}$ and a regular conditional probability $\{(S_y,\mfA),\;y\in{Y}\} $ on $\mfA$ with respect to $Q$ be given. Moreover, let  $\{\xi_y:y\in{Y}\}$ and  $\{\theta_y:y\in{Y}\}$ be stochastic processes defined on $(X,\mfA,P)$. We say that the processes are equivalent if $S_y\{x\in{X}:\xi_y(x)\neq\theta_y(x)\}=0$ for every $y\in{Y}$. The process $\{\theta_y:y\in{Y}\}$ is called measurable if the function $(x,y)\longrightarrow\theta_y(x)$ is $\wh{R}$-measurable.

There are known several examples of stochastic processes without equivalent measurable version (cf. \cite{Cohn, Dudley, HJ}).
The following theorem describes the stochastic processes possessing such a measurable modification. It is an essential strengthening  of \cite[Theorem 5.1]{smm} and  \cite[Theorem 6.1]{mms1}, where only measurable processes were modified.
\begin{thm}
Let $\{\xi_y:y\in{Y}\}$ be a stochastic process defined on $(X,\mfA,P)$. For each $y\in{Y}$ let $S_y$ be the distribution of $\xi_y$ and $R\in{P}\cd{Q}$ be determined by the regular conditional probability $\{S_y:y\in{Y}\}$. Assume one of the following conditions:
\begin{itemize}
\item [$(\alpha)$]
$(X,\mfA,P)$ is separable in the Frechet-Nikod\'{y}m pseudometric;
\item [$(\beta)$]
$R\ll{P}\otimes{Q}$.
\end{itemize}
Then,  the process $\{\xi_y:y\in{Y}\}$  possesses an equivalent measurable version if and only if  the process is $\wh{R_{\dd}}$-measurable. Every modification of a measurable process with the help of  a family of liftings generated by an equi-admissibly generated collection of densities is a measurable process.
\end{thm}
\begin{proof}
Assume that $[\alpha]$ is fulfilled.  According to \cite[Theorem 3.6]{smm} and Lemma \ref{L4} there exist a lifting $\pi:\mfA\dd\mfB\to \mfA\wh\otimes\mfB$ and liftings $\{\sigma_y\in\vL(\wh{S_y}):y\in{Y}\}$ such that $[\pi(f)]^y = \sigma_y\bigl([\pi(f)]^y\bigr)$ for every $f \in {\mcL}^{\infty}(\wh{R_{\dd}})$  and every $y\in Y$.
The same follows from Theorem \ref{T5} in case of $[\beta]$.
The rest of the proof is a copy of the proof of \cite[Proposition 4.1]{mu25}.
\end{proof}
\noindent The full version of \cite[Proposition 4.1]{mu25}  remains a hypothesis.

\end{document}